\documentclass[12pt]{article}
\usepackage{amsmath,amsfonts,amsthm,amssymb}
\newcommand{\F}{\mathbb F_q}
\newcommand{\K}{\overline{K}_c}
\newcommand{\KT}{\overline{K}_c\{ t\}}

\begin{document}
\newtheorem{lem}{Lemma}
\newtheorem{teo}{Theorem}
\newtheorem{def1}{Definition}
\pagestyle{plain}
\title{Umbral Calculus in Positive Characteristic}
\author{Anatoly N. Kochubei\footnote{Partially supported by
CRDF under Grant UM1-2567-OD-03}\\ 
\footnotesize Institute of Mathematics,\\ 
\footnotesize National Academy of Sciences of Ukraine,\\ 
\footnotesize Tereshchenkivska 3, Kiev, 01601 Ukraine
\\ \footnotesize E-mail: \ kochubei@i.com.ua}
\date{}
\maketitle
\newpage
\begin{abstract}
An umbral calculus over local fields of positive characteristic
is developed on the basis of a relation of binomial type 
satisfied by the Carlitz polynomials. Orthonormal bases in the 
space of continuous $\F$-linear functions are constructed.
\end{abstract}
\vspace{2cm}
\newpage

\section{INTRODUCTION}

Classical umbral calculus \cite{RKO,RR,RT} is a set of algebraic 
tools for obtaining, in a unified way, a rich variety of results 
regarding structure and properties of various polynomial 
sequences. There exists a lot of generalizations extending umbral 
methods to other classes of functions. However there is a 
restriction common to the whole literature on umbral calculus -- 
the underlying field must be of zero characteristic. An attempt to 
mimic the characteristic zero procedures in the positive 
characteristic case \cite{F} revealed a number of pathological 
properties of the resulting structures. More importantly, these 
structures were not connected with the existing analysis in 
positive characteristic based on a completely different algebraic 
foundation.

It is well known that any non-discrete locally compact topological 
field of a positive characteristic $p$ is isomorphic to the field 
$K$ of formal Laurent series with coefficients from the Galois 
field $\F$, $q=p^\nu$, $\nu \in \mathbb Z_+$. Denote by $|\cdot |$ 
the non-Archimedean absolute value on $K$; if $z\in K$,
$$
z=\sum\limits_{i=m}^\infty \zeta_ix^i,\quad n\in \mathbb 
Z,\zeta_i\in \F ,\zeta_m\ne 0,
$$
then $|z|=q^{-m}$. This valuation can be extended onto the field 
$\K$, the completion of an algebraic closure of $K$. Let $O=\{ 
z\in K:\ |z|\le 1\}$ be the ring of integers in $K$. The ring $\F 
[x]$ of polynomials (in the indeterminate $x$) with coefficients 
from $\F$ is dense in $O$ with respect to the topology defined by 
the metric $d(z_1,z_2)=|z_1-z_2|$.

It is obvious that standard notions of analysis do not make sense 
in the characteristic $p$ case. For instance, $n!=0$ if $n\ge p$, 
so that one cannot define a usual exponential function on $K$, and 
$\dfrac{d}{dt}(t^n)=0$ if $p$ divides $n$. On the other hand, some 
well-defined functions have unusual properties. In particular, 
there are many functions with the $\F$-linearity property
$$
f(t_1+t_2)=f(t_1)+f(t_2),\quad f(\alpha t)=\alpha f(t),
$$
for any $t_1,t_2,t\in K$, $\alpha \in \F$. Such are, for example, all 
power series $\sum c_kt^{q^k}$, $c_k\in \K$, convergent on some region in
$K$ or $\K$.

The analysis on $K$ taking into account the above special features 
was initiated in a seminal paper by Carlitz \cite{C1} who introduced, for this
case, the appropriate notions of a factorial, an exponential and a logarithm,
a system of polynomials $\{ e_i\}$ (now called the Carlitz polynomials),
and other related objects. In subsequent works by Carlitz, Goss, Thakur, and
many others (see references in \cite{G}) analogs of the gamma, zeta,
Bessel and hypergeometric functions were introduced and studied. A 
difference operator $\Delta$ acting on functions over $K$ or its 
subsets, which was mentioned briefly in \cite{C1}, became (as an 
analog of the operator $t\frac{d}{dt}$) the main ingredient of the 
calculus and the analytic theory of differential equations on $K$ 
\cite{K2,K3,K4}. It appears also in a characteristic $p$ analog of 
the canonical commutation relations of quantum mechanics found in 
\cite{K1}.

The definition of the Carlitz polynomials is as follows. Let 
$e_0(t)=t$,
\begin{equation}
e_i(t)=\prod \limits _{\genfrac{}{}{0pt}{1}{m\in \F [x]}{\deg m<i}}(t-m),\quad
i\ge 1
\end{equation}
(we follow the notation in \cite{G0} used in the modern 
literature; the initial formulas from \cite{C1} have different 
signs in some places). It is known \cite{C1,G0} that 
\begin{equation}
e_i(t)=\sum \limits _{j=0}^i(-1)^{i-j}\genfrac{[}{]}{0pt}{0}{i}{j}
t^{q^j}
\end{equation}
where
$$
\genfrac{[}{]}{0pt}{0}{i}{j}=\frac{D_i}{D_jL_{i-j}^{q^j}},
$$
$D_i$ is the Carlitz factorial
\begin{equation}
D_i=[i][i-1]^q\ldots [1]^{q^{i-1}},\quad [i]=x^{q^i}-x\ (i\ge 1),\ 
D_0=1,
\end{equation}
the sequence $\{L_i\}$ is defined by
\begin{equation}
L_i=[i][i-1]\ldots [1]\ (i\ge 1);\ L_0=1.
\end{equation}
It follows from (3), (4) that
$$
|D_i|=q^{-\frac{q^i-1}{q-1}},\quad |L_i|=q^{-i}.
$$

The normalized polynomials $f_i(t)=\dfrac{e_i(t)}{D_i}$ form an 
orthonormal basis in the Banach space $C_0(O,\K )$ of all 
$\F$-linear continuous functions $O\to \K$, with the supremum norm $\| \cdot 
\|$. Thus every function $\varphi \in C_0(O,\K )$ admits a unique 
representation as a uniformly convergent series
$$
\varphi =\sum\limits_{i=0}^\infty a_if_i,\quad a_i\in \K ,\ a_i\to 
0,
$$
satisfying the orthonormality condition
$$
\varphi =\sup\limits_{i\ge 0}|a_i|.
$$
For several different proofs of this fact see \cite{Con,K1,W}. 
Note that we consider functions with values in $\K$ defined on a 
compact subset of $K$.

The sequences $\{ D_i\}$ and $\{ L_i\}$ are involved in the 
definitions of the Carlitz exponential and logarithm:
\begin{equation}
e_C(t)=\sum\limits_{n=0}^\infty \frac{t^{q^n}}{D_n},\quad 
\log_C(t)=\sum\limits_{n=0}^\infty (-1)^n\frac{t^{q^n}}{L_n},\quad 
|t|<1.
\end{equation}
It is seen from (2) and (5) that the above functions are 
$\F$-linear on their domains of definition.

The most important object connected with the Carlitz polynomials 
is the Carlitz module
\begin{equation}
C_s(z)=\sum\limits_{i=0}^{\deg s}f_i(s)z^{q^i}=\sum\limits_{i=0}^{\deg 
s}\frac{e_i(s)}{D_i}z^{q^i},\quad s\in \F [x].
\end{equation}
Note that by (1) $e_i(s)=0$ if $s\in \F [x]$, $\deg s<i$.

The function $C_s$ appears in the functional equation for the 
Carlitz exponential,
$$
C_s(e_C(t))=e_C(st).
$$
Its main property is the relation
\begin{equation}
C_{ts}(z)=C_t(C_s(z)),\quad s,t\in \F [x],
\end{equation}
which obtained a far-reaching generalization in the theory of 
Drinfeld modules, the principal objects of the function field 
arithmetic (see \cite{G}).

Let us write the identity (7) explicitly using (6). After 
rearranging the sums we find that
$$
C_{ts}(z)=\sum\limits_{i=0}^{\deg t+\deg s}z^{q^i}\sum
\limits_{\genfrac{}{}{0pt}{1}{m+n=l}{m,n\ge 
0}}\frac{1}{D_nD_m^{q^n}}e_n(t)\{ e_m(s)\}^{q^n},
$$
so that 
\begin{equation}
e_i(st)=\sum\limits_{n=0}^i\binom{i}{n}_Ke_n(t)\{ e_{i-n}(s)\}^{q^n}
\end{equation}
where
\begin{equation}
\binom{i}{n}_K=\frac{D_i}{D_nD_{i-n}^{q^n}}.
\end{equation}

In this paper we show that the ``$K$-binomial'' relation (8), a 
positive characteristic counterpart of the classical binomial 
formula, can be used for developing umbral calculus in the spirit 
of \cite{RKO}. In particular, we introduce and study corresponding 
(nonlinear) delta operators, obtain a representation for operators 
invariant with respect to multiplicative shifts, construct generating 
functions for polynomial sequences of the $K$-binomial type. Such 
sequences are also used for constructing new orthonormal bases of 
the space $C_0(O,\K )$ (in particular, a sequence of the Laguerre 
type polynomials), in a way similar to the $p$-adic 
(characteristic 0) case \cite{VH,Ver,Rob}.

\section{DELTA OPERATORS AND $K$-BINOMIAL \\ SEQUENCES}

Denote by $\KT$ the vector space over $\K$ consisting of 
$\F$-linear polynomials $u=\sum a_kt^{q^k}$ with coefficients from 
$\K$. We will often use the operator of multiplicative shift 
$(\rho_\lambda u)(t)=u(\lambda t)$ on $\KT$ and the Frobenius 
operator $\tau u=u^q$. We call a linear operator $T$ on $\KT$ {\it 
invariant} if it commutes with $\rho_\lambda$ for each $\lambda 
\in K$.

\medskip
\begin{lem}
If $T$ is an invariant operator, then $T(t^{q^n})=c_nt^{q^n}$, 
$c_n\in \K$, for each $n\ge 0$.
\end{lem}

\medskip
{\it Proof}. Suppose that
$$
T(t^{q^n})=\sum\limits_{l=1}^Nc_{j_l}t^{q^{j_l}}
$$
where $j_l$ are different non-negative integers, $c_{j_l}\in \K$. 
For any $\lambda \in K$
$$
\rho_\lambda T(t^{q^n})=T\rho_\lambda (t^{q^n})=T\left( (\lambda 
t)^{q^n}\right) 
=\lambda^{q^n}T(t^{q^n})=\lambda^{q^n}\sum\limits_{l=1}^Nc_{j_l}t^{q^{j_l}}.
$$

On the other hand,
$$
\rho_\lambda 
T(t^{q^n})=\sum\limits_{l=1}^Nc_{j_l}\lambda^{q^{j_l}}t^{q^{j_l}}.
$$
Since $\lambda$ is arbitrary, this implies the required result. 
$\blacksquare$

\medskip
If an invariant operator $T$ is such that $T(t)=0$, then by Lemma 
1 the operator $\tau^{-1}T$ on $\KT$ is well-defined.

\medskip
\begin{def1}
A $\F$-linear operator $\delta =\tau^{-1}\delta_0$, where 
$\delta_0$ is a linear invariant operator on $\KT$, is called a 
delta operator if $\delta_0(t)=0$ and $\delta_0(f)\ne 0$ for $\deg 
f>1$, that is $\delta_0(t^{q^n})=c_nt^{q^n}$, $c_n\ne 0$, for all 
$n\ge 1$.
\end{def1}

\medskip
The most important example of a delta operator is {\it the Carlitz 
derivative} $d=\tau^{-1}\Delta$ where
$$
(\Delta u)(t)=u(xt)-xu(t).
$$
Many interesting $\F$-linear functions satisfy equations involving 
the operator $d$; for example, for the Carlitz exponential we have 
$de_C=e_C$. It appears also in $\F$-linear representations of the 
canonical commutation relations \cite{K1,K2}.

\medskip
\begin{def1}
A sequence $\{ P_n\}_0^\infty$ of $\F$-linear polynomials is 
called a basic sequence corresponding to a delta operator $\delta 
=\tau^{-1}\delta_0$, if $\deg P_n=q^n$, $P_0(1)=1$, $P_n(1)=0$ for 
$n\ge 1$,
\begin{equation}
\delta P_0=0,\quad \delta P_n=[n]^{1/q}P_{n-1},\ n\ge 1,
\end{equation}
or, equivalently,
\begin{equation}
\delta_0P_0=0,\quad \delta_0P_n=[n]P_{n-1}^q,\ n\ge 1.
\end{equation}
\end{def1}

\medskip
It follows from well-known identities for the Carlitz polynomials 
$e_i$ (see \cite{G0}) that the sequence $\{ e_i\}$ is basic with 
respect to the operator $d$. For the normalized Carlitz polynomials 
$f_i$ we have the relations 
$$
df_0=0,\quad df_i=f_{i-1},\ i\ge 1.
$$ 

The next definition is a formalization of the property (8).

\medskip
\begin{def1}
A sequence of $\F$-linear polynomials $u_i\in \KT$ is called a 
sequence of $K$-binomial type if $\deg u_i=q^i$ and for all 
$i=0,1,2,\ldots$
\begin{equation}
u_i(st)=\sum\limits_{n=0}^i\binom{i}{n}_Ku_n(t)\left\{ 
u_{i-n}(s)\right\}^{q^n},\quad s,t\in K.
\end{equation}
\end{def1}

\medskip
If $\{ u_i\}$ is a sequence of $K$-binomial type, then $u_i(1)=0$ 
for $i\ge 1$, $u_0(1)=1$ (so that $u_0(t)=t$).

Indeed, for $i=0$ the formula (12) gives $u_0(st)=u_0(s)u_0(t)$. 
Setting $s=1$ we have $u_0(t)=u_0(1)u_0(t)$, and since $\deg 
u_0=1$, so that $u_0(t)\not \equiv 0$, we get $u_0(1)=1$.

If $i>0$, for all $t$
$$
0=u_i(t)-u_i(t)=\sum\limits_{n=0}^{i-1}\binom{i}{n}_K\left\{ 
u_{i-n}(1)\right\}^{q^n}u_n(t),
$$
and the linear independence of the polynomials $u_n$ means that 
$u_l(1)=0$ for $l\ge 1$.

\medskip
\begin{teo}
For any delta operator $\delta =\tau^{-1}\delta_0$, there exists a 
unique basic sequence $\{ P_n\}$, which is a sequence of $K$-binomial 
type. Conversely, given a sequence $\{ P_n\}$ of $K$-binomial 
type, define the action of $\delta_0$ on $P_n$ by the relations (11), 
extend it onto $\KT$ by linearity and set $\delta =\tau^{-1}\delta_0$. Then
$\delta$ is a delta operator, and $\{ P_n\}$ is the corresponding 
basic sequence.
\end{teo}

{\it Proof}. Let us construct a basic sequence corresponding to 
$\delta$. Set $P_0(t)=t$ and suppose that $P_{n-1}$ has been 
constructed. It follows from Lemma 1 that $\delta$ is surjective, 
and we can choose $P_n$ satisfying (10). For any $c\in \K$, 
$P_n+ct$ also satisfies (10), and we may redefine $P_n$ choosing 
$c$ in such a way that $P_n(1)=0$.

Hence, a basic sequence $\{ P_n\}$ indeed exists. If there is 
another basic sequence $\{ P_n'\}$ with the same delta operator, 
then $\delta (P_n-P_n')=0$, whence $P_n'(t)=P_n(t)+at$, $a\in \K$, 
and setting $t=1$ we find that $a=0$.

In order to prove the $K$-binomial property, we introduce some 
operators having an independent interest.

Consider the linear operators $\delta_0^{(l)}=\tau^l\delta^l$.

\medskip
\begin{lem}
\begin{description}
\item[(i)] The identity
\begin{equation}
\delta_0^{(l)}P_j=\frac{D_j}{D_{j-l}^{q^l}}P_{j-l}^{q^l}
\end{equation}
holds for any $l\le j$.
\item[(ii)] Let $f$ be a $\F$-linear polynomial, $\deg f\le q^n$. 
Then a generalized Taylor formula
\begin{equation}
f(st)=\sum\limits_{l=0}^n\frac{\left( \delta_0^{(l)}f\right) 
(s)}{D_l}P_l(t)
\end{equation}
holds for any $s,t\in K$.
\end{description}
\end{lem}

{\it Proof}. By (10),
\begin{multline*}
\delta^lP_j=\delta^{l-1}\left( [j]^{q^{-1}}P_{j-1}\right) 
=[j]^{q^{-l}}\delta^{l-1}P_{j-1}
=[j]^{q^{-l}}[j-1]^{q^{-(l-1)}}\delta^{l-2}P_{j-2}=\ldots \\
=[j]^{q^{-l}}[j-1]^{q^{-(l-1)}}\ldots [j-(l-1)]^{q^{-1}}P_{j-l},
\end{multline*}
so that
$$
\delta_0^{(l)}P_j=[j][j-1]^q\ldots [j-(l-1)]^{q^{l-1}}P_{j-l}^{q^l}
$$
which is equivalent to (13).

Since $\deg P_j=q^j$, the polynomials $P_1,\ldots ,P_n$ form a basis 
of the vector space of all $\F$-linear polynomials of degrees $\le 
n$ (because its dimension equals $n$). Therefore
\begin{equation}
f(st)=\sum\limits_{j=0}^nb_j(s)P_j(t)
\end{equation}
where $b_j(s)$ are, for each fixed $s$, some elements of $\K$.

Applying the operator $\delta_0^{(l)}$, $0\le l\le n$, in the 
variable $t$ to both sides of (15) and using (14) we find that
$$
\left( \delta_0^{(l)}f\right) 
(st)=\sum\limits_{j=l}^nb_j(s)\frac{D_j}{D_{j-l}^{q^l}}P_{j-l}^{q^l}(t)
$$
(note also that $\delta_0^{(l)}$ commutes with $\rho_s$). Setting $t=1$ 
and taking into account that
$$
P_{j-l}(1)=\begin{cases}
0, & \text{if $j>l$;}\\
1, & \text{if $j=l$,}
\end{cases}
$$
we come to the equality
$$
b_l(s)=\frac{\left( \delta_0^{(l)}f\right) (s)}{D_l},\quad 0\le 
l\le n,
$$
which implies (14). $\blacksquare$

\medskip
Note that the formulas (13) and (14) for the Carlitz polynomials 
$e_i$ are well known; see \cite{G0}. It is important that, in 
contrast to the classical umbral calculus, the linear operators 
involved in (14) are not powers of a single linear operator.

\medskip
{\it Proof of Theorem 1 (continued)}. In order to prove that $\{ 
P_n\}$ is a sequence of $K$-binomial type, it suffices to take 
$f=P_n$ in (14) and to use the identity (13).

To prove the second part of the theorem, we calculate the action 
in the variable $t$ of the operator $\delta_0$, defined by (11), 
upon the function $P_n(st)$. Using the relation 
$D_{n+1}=[n+1]D_n^q$ we find that
\begin{multline*}
\delta_{0,t}P_n(st)=\sum\limits_{j=0}^n\binom{n}{j}_KP_{n-j}^{q^j}
(s)\left( \delta_0P_j\right) 
(t)=\sum\limits_{j=1}^n\frac{D_n}{D_{n-j}^{q^j}D_j}P_{n-j}^{q^j}
(s)[j]P_{j-1}^q(t)\\
=\sum\limits_{i=0}^{n-1}\frac{D_n[i+1]}{D_{n-i-1}^{q^{i+1}}
D_{i+1}}P_{n-i-1}^{q^{i+1}}(s)P_i^q(t)=[n]\sum\limits_{i=0}^{n-1}
\left( 
\frac{D_{n-1}}{D_{n-i-1}^qD_i}\right)^qP_{n-i-1}^{q^{i+1}}(s)P_i^q(t)\\
=[n]\left\{ 
\sum\limits_{i=0}^{n-1}\binom{n-1}{i}_KP_{n-i-1}^{q^i}(s)P_i(t)\right\}^q
=[n]P_{n-1}^q(st)=\left( \delta_0P_n\right) (st),
\end{multline*}
that is $\delta_0$ commutes with multiplicative shifts.

It remains to prove that $\delta_0(f)\ne 0$ if $\deg f>1$. 
Assuming that $\delta_0(f)=0$ for $f=\sum\limits_{j=0}^na_jP_j$ we 
have
$$
0=\sum\limits_{j=0}^na_j[j]P_{j-1}^q=\left\{ 
\sum\limits_{i=0}^{n-1}a_{i+1}^{1/q}[i+1]^{1/q}P_i\right\}^q
$$
whence $a_1=a_2=\ldots =a_n=0$ due to the linear independence of 
the sequence $\{ P_i\}$. $\blacksquare$

\section{INVARIANT OPERATORS}

Let $T$ be a linear invariant operator on $\KT$. Let us find its 
representation via an arbitrary fixed delta operator $\delta 
=\tau^{-1}\delta_0$. By (14), for any $f\in \KT$, $\deg f=q^n$,
$$
(Tf)(st)=\left( \rho_sTf\right) (t)=T\left( \rho_sf\right) 
(t)=T_tf(st)=\sum\limits_{l=0}^n\left( TP_l\right) (t)\frac{\left( 
\delta_0^{(l)}f\right) (s)}{D_l}.
$$
Setting $s=1$ we find that
\begin{equation}
T=\sum\limits_{l=0}^\infty \sigma_l\delta_0^{(l)}
\end{equation}
where $\sigma_l=\dfrac{(TP_l)(1)}{D_l}$. The infinite series in 
(16) becomes actually a finite sum if both sides of (16) are 
applied to any $\F$-linear polynomial $f\in \KT$.

Conversely, any such series defines a linear invariant operator on 
$\KT$.

Below we will consider in detail the case where $\delta$ is the 
Carlitz derivative $d$, so that $\delta_0=\Delta$, and the 
operators $\delta_0^{(l)}=\Delta^{(l)}$ are given recursively 
\cite{G0}:
\begin{equation}
\left( \Delta^{(l)}u\right) (t)=\left( \Delta^{(l-1)}u\right) 
(xt)-x^{q^{l-1}}u(t)
\end{equation}
(the formula (16) for this case was proved by a different method 
in \cite{J}).

Using (17) with $l=0$, we can compute for this case the 
coefficients $c_n$ from Lemma 1. We have 
$\Delta^{(l)}(t^{q^n})=0$, if $n<l$,
$$
\Delta (t^{q^n})=[n]t^{q^n},\quad n\ge 1;
$$
$$
\Delta^{(2)}(t^{q^n})=\tau^2d^2(t^{q^n})=\tau \Delta \tau^{-1}\Delta 
(t^{q^n})=\tau \Delta \left( [n]^{1/q}t^{q^{n-1}}\right) 
=[n][n-1]^qt^{q^n},\quad n\ge 2,
$$
and by induction
\begin{equation}
\Delta^{(l)}(t^{q^n})[n][n-1]^q\ldots 
[n-l+1]^{q^{l-1}}t^{q^n}=\frac{D_n}{D_{n-l}^{q^l}}t^{q^n},\quad n\ge l.
\end{equation}

The explicit formula (18) makes it possible to find out when an 
operator $\theta =\tau^{-1}\theta_0$, with
\begin{equation}
\theta_0=\sum\limits_{l=1}^\infty \sigma_l\Delta^{(l)},
\end{equation}
is a delta operator. We have $\theta_0(t)=0$,
$$
\theta_0(t^{q^n})=D_nS_nt^{q^n},
$$
where $S_n=\sum\limits_{l=1}^n\dfrac{\sigma_l}{D_{n-l}^{q^l}}$. 
Thus $\theta$ is a delta operator if and only if $S_n\ne 0$ for 
all $n=1,2,\ldots$.

\medskip
{\bf Example 1}.
Let $\sigma_l=1$ for all $l\ge 1$, that is
\begin{equation}
\theta_0=\sum\limits_{l=1}^\infty \Delta^{(l)}.
\end{equation}

Since $|D_i|=q^{-\frac{q^i-1}{q-1}}$, we have
$$
\left| D_{n-l}^{q^l}\right| =q^{-\frac{q^n-q^l}{q-1}},
$$
so that $|S_n|=q^{\frac{q^n-q}{q-1}}$ ($\ne 0$) by the 
ultra-metric property of the absolute value. Comparing (20) with a 
classical formula from \cite{RKO} we may see the polynomials $P_n$ 
for this case as analogs of the Laguerre polynomials.

{\bf Example 2}. Let $\sigma_l=\dfrac{(-1)^{l+1}}{L_l}$. Now
$$
S_n=\sum\limits_{l=1}^n(-1)^{l+1}\dfrac{1}{L_lD_{n-l}^{q^l}}.
$$
Let us use the identity
\begin{equation}
\sum\limits_{j=0}^{h-1}\frac{(-1)^j}{L_jD_{h-j}^{q^j}}=
\frac{(-1)^{h+1}}{L_h}
\end{equation}
proved in \cite{Gek}. It follows from (21) that
$$
\sum\limits_{j=1}^h\frac{(-1)^j}{L_jD_{h-j}^{q^j}}=
\frac{(-1)^{h+1}}{L_h}-\frac{1}{D_h}+\frac{(-1)^h}{L_h}
=-\frac{1}{D_h},
$$
so that $S_n=D_n^{-1}$ ($\ne 0$), $n=1,2,\ldots$. In this case 
$\theta_0(t^{q^j})=t^{q^j}$ for all $j\ge 1$ (of course, 
$\theta_0(t)=0$), and $P_0(t)=t$, $P_n(t)=D_n\left( 
t^{q^n}-t^{q^{n-1}}\right)$ for $n\ge 1$.

\section{ORTHONORMAL BASES}

Let $\{ P_n\}$ be the basic sequence corresponding to a delta 
operator $\delta =\tau^{-1}\delta_0$,
$$
\delta_0=\sum\limits_{l=1}^\infty \sigma_l\Delta^{(l)}
$$
(the operator series converges on any polynomial from $\KT$).

Let $Q_n=\dfrac{P_n}{D_n}$, $n=0,1,2,\ldots$. Then for any $n\ge 
1$
$$
\delta Q_n=D_n^{-1/q}\delta 
P_n=\frac{[n]^{1/q}}{D_n^{1/q}}P_{n-1}=\frac{P_{n-1}}{D_{n-1}}=Q_{n-1},
$$
and the $K$-binomial property of $\{ P_n\}$ implies the identity
\begin{equation}
Q_i(st)=\sum\limits_{n=0}^iQ_n(t)\left\{ 
Q_{i-n}(s)\right\}^{q^n},\quad s,t\in K.
\end{equation}

The identity (22) may be seen as another form of the $K$-binomial 
property. Though it resembles its classical counterpart, the 
presence of the Frobenius powers is a feature specific for the 
case of a positive characteristic. We will call $\{ Q_n\}$ {\it a 
normalized basic sequence}.

\medskip
\begin{teo}
If $|\sigma_1|=1$, $|\sigma_l|\le 1$ for $l\ge 2$, then the 
sequence $\{ Q_n\}_0^\infty$ is an orthonormal basis of the space 
$C_0(O,\K )$ -- for any $f\in C_0(O,\K )$ there is a uniformly 
convergent expansion
\begin{equation}
f(t)=\sum\limits_{n=0}^\infty \psi_nQ_n(t),\quad t\in O,
\end{equation}
where $\psi_n=\left( \delta_0^{(n)}f\right) (1)$, $|\psi_n|\to 0$ 
as $n\to \infty$,
\begin{equation}
\|f\| =\sup\limits_{n\ge 0}|\psi_n|.
\end{equation}
\end{teo}

\medskip
{\it Proof}. We have $Q_0(t)=P_0(t)=t$, so that $\|Q_0\|=1$. Let 
us prove that $\|Q_n\|=1$ for all $n\ge 1$. Our reasoning will be 
based on expansions in the normalized Carlitz polynomials $f_n$.

Let $n=1$. Since $\deg Q_n=q^n$, we have $Q_1=a_0f_0+a_1f_1$. We 
know that $Q_1(1)=f_1(1)=0$, hence $a_0=0$, so that $Q_1=a_1f_1$. 
Next, $\delta Q_1=Q_0=f_0$. Writing this explicitly we find that
$$
f_0=a_1^{1/q}\tau^{-1}\sum\limits_{l=1}^\infty 
\sigma_l\Delta^{(l)}f_1=a_1^{1/q}\sigma_1^{1/q}df_1=a_1^{1/q}\sigma_1^{1/q}f_0,
$$
whence $a_1=\sigma_1^{-1}$, $Q_1=\sigma_1^{-1}f_1$, and 
$\|Q_1\|=1$.

Assume that $\|Q_{n-1}\|=1$ and consider the expansion
$$
Q_n=\sum\limits_{j=1}^na_jf_j
$$
(the term containing $f_0$ is absent since $Q_n(1)=0$). Applying 
$\delta$ we get
$$
\delta Q_n=\sum\limits_{j=1}^na_j^{1/q}\sum\limits_{l=1}^\infty 
\sigma_l^{1/q}\tau^{-1}\Delta^{(l)}f_j.
$$

It is known \cite{G0} that
$$
\Delta^{(l)}e_j=\begin{cases}
\frac{D_j}{D_{j-l}^{q^l}}{e_{j-l}^{q^l}}, & \text{if $l\le j$},\\
0, & \text{if $l>j$},\end{cases}
$$
so that
$$
\Delta^{(l)}f_j=\begin{cases}
f_{j-l}^{q^l}, & \text{if $l\le j$},\\
0, & \text{if $l>j$}.\end{cases}
$$
Therefore
\begin{equation}
\delta Q_n=\sum\limits_{j=1}^na_j^{1/q}\sum\limits_{l=1}^j 
\sigma_l^{1/q}f_{j-l}^{q^{l-1}}.
\end{equation}

It follows from the identity $f_{i-1}^q=f_{i-1}+[i]f_i$ (see 
\cite{G0,K1}) that
$$
f_{j-l}^{q^{l-1}}=\sum\limits_{k=0}^{l-1}\varphi_{j,l,k}f_{j-l+k}
$$
where $\varphi_{j,l,0}=1$, $\left| \varphi_{j,l,k}\right| <1$ for 
$k\ge 1$. Substituting into (25) we find that
\begin{multline*}
Q_{n-1}=\sum\limits_{j=1}^na_j^{1/q}\sum\limits_{l=1}^j 
\sigma_l^{1/q}\sum\limits_{k=0}^{l-1}\varphi_{j,l,k}f_{j-l+k}
=\sum\limits_{j=1}^na_j^{1/q}\sum\limits_{i=0}^{j-1}f_i
\sum\limits_{l=j-i}^j\sigma_l^{1/q}\varphi_{j,l,i-j+l}\\
=\sum\limits_{i=0}^{n-1}f_i\sum\limits_{j=i+1}^na_j^{1/q}
\sum\limits_{l=j-i}^j\sigma_l^{1/q}\varphi_{j,l,i-j+l}
\end{multline*}
whence
\begin{equation}
\max\limits_{0\le i\le n-1}\left| \sum\limits_{j=i+1}^na_j^{1/q}
\sum\limits_{l=j-i}^j\sigma_l^{1/q}\varphi_{j,l,i-j+l}\right| =1
\end{equation}
by the inductive assumption and the orthonormal basis property of 
the normalized Carlitz polynomials.

For $i=n-1$, we obtain from (26) that
$$
\left| 
a_n^{1/q}\sum\limits_{l=1}^n\sigma_l^{1/q}\varphi_{n,l,l-1}\right| 
\le 1.
$$
We have $\varphi_{n,1,0}=1$, $|\sigma_1|=1$, and
$$
\left| \sum\limits_{l=2}^n\sigma_l^{1/q}\varphi_{n,l,l-1}\right| 
<1,
$$
so that
$$
\left| \sum\limits_{l=1}^n\sigma_l^{1/q}\varphi_{n,l,l-1}\right| 
=1
$$
whence $|a_n|\le 1$.

Next, for $i=n-2$ we find from (26) that
$$
\left| 
a_{n-1}^{1/q}\sum\limits_{l=1}^{n-1}\sigma_l^{1/q}\varphi_{n-1,l,l-1}+
a_n^{1/q}\sum\limits_{l=2}^n\sigma_l^{1/q}\varphi_{n,l,l-1}\right| 
\le 1.
$$
We have proved that the second summand on the left is in $O$; then 
the first summand is considered as above, so that $|a_{n-1}|\le 
1$. Repeating this reasoning we come to the conclusion that 
$|a_j|\le 1$ for all $j$. Moreover, $|a_j|=1$ for at least one 
value of $j$; otherwise we would come to a contradiction with 
(26). This means that $\| Q_n\| =1$.

If $f$ is an arbitrary $\F$-linear polynomial, $\deg f=q^N$, then 
by the generalized Taylor formula (14)
$$
f(t)=\sum\limits_{l=0}^N\psi_lQ_l(t),\quad t\in O,
$$
where $\psi_l=\left( \delta_0^{(l)}f\right) (1)$.

Since $\|Q_l\|=1$ for all $l$, we have $\|f\|\le 
\sup\limits_l|\psi_l|$. On the other hand, 
$\delta_0^{(l)}f=\tau^l\left( \tau^{-1}\delta_0\right)^lf$, and 
if we prove that $\| \delta_0f\| \le \| f\|$, this will imply the 
inequality $\| \delta_0^{(l)}f\| \le \| f\|$. We have
\begin{multline*}
\left\| \Delta^{(l)}f\right\| =\max\limits_{t\in O}\left| \left(
\Delta^{(l-1)}f\right) (xt)-x^{q^{l-1}}\left(
\Delta^{(l-1)}f\right) (t)\right| \le \max\limits_{t\in O}\left| \left(
\Delta^{(l-1)}f\right) (t)\right| \le \ldots \\
\le \max\limits_{t\in O}\left| \left( \Delta f\right) (t)\right| 
\le \|f\|,
\end{multline*}
so that
$$
\| \delta_0f\|=\left\| \sum\limits_{l=0}^\infty 
\sigma_l\Delta^{(l)}f\right\| \le \sup\limits_l|\sigma_l|\cdot
\left\| \Delta^{(l)}f\right\| \le \|f\|
$$
whence $\| \delta_0^{(l)}f\| \le \| f\|$ and 
$\sup\limits_l|\psi_l|\le \|f\|$.

Thus, we have proved (24) for any polynomial. By a well-known 
result of non-Archimedean functional analysis (see Theorem 50.7 in 
\cite{Sch}), the uniformly convergent expansion (23) and the 
equality (24) hold for any $f\in C_0(O,\K )$. 

The relation $\psi_n=\left( \delta_0^{(n)}f\right) (1)$ also 
remains valid for any $f\in C_0(O,\K )$. Indeed, denote by 
$\varphi_n(f)$ a continuous linear functional on $C_0(O,\K )$ of 
the form $\left( \delta_0^{(n)}f\right) (1)$. Suppose that 
$\{F_N\}$ is a sequence of $\F$-linear polynomials uniformly 
convergent to $f$. Then
$$
F_N=\sum\limits_n\varphi_n(F_N)Q_n,
$$
so that
$$
F-F_N=\sum\limits_{n=0}^\infty \left\{ 
\psi_n-\varphi_n(F_N)\right\}Q_n,
$$
and by (24),
$$
\left\| F-F_N\right\| =\sup\limits_n \left| 
\psi_n-\varphi_n(F_N)\right| .
$$

For each fixed $n$ we find that $\left| 
\psi_n-\varphi_n(F_N)\right| \le \left\|F-F_N\right\|$, and 
passing to the limit as $N\to \infty$ we get that 
$\psi_n=\varphi_n(f)$, as desired. $\blacksquare$

\medskip
By Theorem 2, the Laguerre-type polynomial sequence from 
Example 1 is an orthonormal basis of $C_0(O,\K )$. The sequence 
from Example 2 does not satisfy the conditions of Theorem 2.

Note that the conditions $|\sigma_1|=1$, $|\sigma_l|\le 1$, 
$l=2,3,\ldots$, imply that $S_n\ne 0$ for all $n$, so that the 
series (19) considered in Theorem 2 always correspond to delta 
operators.

Let us write a recurrence formula for the coefficients of the 
polynomials $Q_n$. Here we assume only that $S_n\ne 0$ for all 
$n$. Let
\begin{equation}
Q_n(t)=\sum\limits_{j=0}^n\gamma_j^{(n)}t^{q^j}.
\end{equation}
We know that $\gamma_0^{(0)}=1$.

Using the relation $\delta_0\left( t^{q^n}\right) =D_nS_nt^{q^n}$ 
we find that for $n\ge 1$
$$
Q_{n-1}=\delta 
Q_n=\tau^{-1}\sum\limits_{j=1}^n\gamma_j^{(n)}D_jS_jt^{q^j}=
\sum\limits_{i=0}^{n-1}\left( 
\gamma_{i+1}^{(n)}\right)^{1/q}D_{i+1}^{1/q}S_{i+1}^{1/q}t^{q^i}.
$$
Comparing this with the equality (27), with $n-1$ substituted for 
$n$, we get
$$
\gamma_i^{(n-1)}=\left( 
\gamma_{i+1}^{(n)}\right)^{1/q}D_{i+1}^{1/q}S_{i+1}^{1/q}
$$
whence
\begin{equation}
\gamma_{i+1}^{(n)}=\frac{\left( 
\gamma_i^{(n-1)}\right)^q}{D_{i+1}S_{i+1}},\quad i=0,1,\ldots 
,n-1;\ n=1,2,\ldots .
\end{equation}

The recurrence formula (28) determines all the coefficients 
$\gamma_i^{(n)}$ (if the polynomial $Q_{n-1}$ is already known) 
except $\gamma_0^{(n)}$. The latter can be found from the 
condition $Q_n(1)=0$:
$$
\gamma_0^{(n)}=-\sum\limits_{j=1}^n\gamma_j^{(n)}.
$$

\section{GENERATING FUNCTIONS}

The definition (6) of the Carlitz module can be seen as a 
generating function for the normalized Carlitz polynomials $f_i$. 
Here we give a similar construction for the normalized basic 
sequence in the general case. As in Sect. 4, we consider a delta 
operator of the form $\delta =\tau^{-1}\delta_0$,
$$
\delta_0=\sum\limits_{l=1}^\infty \sigma_l\Delta^{(l)}
$$
We assume that $S_n\ne 0$ for all $n$.

Let us define the generalized exponential
\begin{equation}
e_\delta (t)=\sum\limits_{j=0}^\infty b_jt^{q^j}
\end{equation}
by the conditions $\delta e_\delta =e_\delta$, $b_0=1$. 
Substituting (29) we come to the recurrence relation
\begin{equation}
b_{j+1}=\frac{b_j^q}{D_{j+1}S_{j+1}}
\end{equation}
which determines $e_\delta$ as a formal power series.

Since $b_0=1$, the composition inverse $\log_\delta$ to the formal power 
series $e_\delta$ has a similar form:
\begin{equation}
\log_\delta (t)=\sum\limits_{n=0}^\infty \beta_nt^{q^n},\quad 
\beta_n\in K,
\end{equation}
(see Sect. 19.7 in \cite{P} for a general treatment of formal 
power series of this kind). A formal substitution gives the 
relations
$$
\beta_0=1,\quad \sum\limits_{m+n=l}b_m\beta_n^{q^m}=0,\quad 
l=1,2,\ldots ,
$$
whence
\begin{equation}
\beta_l=-\sum\limits_{m=1}^lb_m\beta_{l-m}^{q^m},\quad l=1,2,\ldots .
\end{equation}

\medskip
\begin{teo}
Suppose that $|\sigma_1|=1$ and $|\sigma_l|\le 1$ for all $l$. 
Then both the series (29) and (31) converge on the disk 
$D_q=\left\{ t\in O:\ |t|\le q^{-1}\right\}$, if $q\ne 2$, or
$D_2=\left\{ t\in O:\ |t|\le q^{-2}\right\}$, if $q=2$, and 
\begin{equation}
e_\delta (t\log_\delta z)=\sum\limits_{n=0}^\infty Q_n(t)z^{q^n},\quad 
t\in O,\ z\in D_q.
\end{equation}
\end{teo}

{\it Proof}. Since
$$
\left| \frac{D_n}{D_{n-l}^{q^l}}\right| =q^{-\frac{q^l-1}{q-1}},
$$
under our assumptions we have $|D_nS_n|=q^{-1}$ for all $n$. By 
(30), $|b_{j+1}|=q|b_j|^q$, $j=0,1,2,\ldots$, and we prove easily by 
induction that
\begin{equation}
|b_j|=q^{\frac{q^j-1}{q-1}},\quad j=0,1,2,\ldots .
\end{equation}

For the sequence (32) we obtain the estimate
\begin{equation}
|\beta_j|\le q^{\frac{q^j-1}{q-1}},\quad j=0,1,2,\ldots .
\end{equation}
Indeed, this is obvious for $j=0$. If (35) is proved for $j\le 
l-1$, then
$$
|\beta_l|\le \max\limits_{1\le m\le l}|b_m|\cdot 
|\beta_{l-m}|^{q^m}\le \max\limits_{1\le m\le 
l}q^{\frac{q^m-1}{q-1}+q^m\frac{q^{l-m}-1}{q-1}}=q^{\frac{q^l-1}{q-1}}.
$$

It follows from (34) and (35) that both the series (29) and (31) 
are convergent for $t\in D_q$ (in fact they are convergent on a 
wider disk from $\K$, but here we consider them only on $K$). Note 
also that
$$
|\log_\delta (t)|\le \max\limits_{n\ge 0}q^{-\frac{1}{q-1}}\left( 
q^{\frac{1}{q-1}}|t|\right)^{q^n} =|t|
$$
if $t\in D_q$.

If $\lambda \in D_q$, then the function $t\mapsto e_\delta 
(\lambda t)$ is continuous on $O$, and by Theorem 2
$$
e_\delta (\lambda t)=\sum\limits_{n=0}^\infty \psi_n(\lambda 
)Q_n(t)
$$
where $\psi_n(\lambda )=\left( \delta_0^{(n)}e_\delta (\lambda 
\cdot )\right) (1)=\left( \delta_0^{(n)}e_\delta \right) (\lambda )
$
due to the invariance of the operator $\delta_0^{(n)}$. Since
$\delta_0^{(n)}=\tau^n\delta^n$ and $\delta e_\delta =e_\delta$, 
we find that $\delta_0^{(n)}e_\delta =e_\delta^{q^n}$. Therefore
\begin{equation}
e_\delta (\lambda t)=\sum\limits_{n=0}^\infty Q_n(t)\left\{ 
e_\delta (\lambda )\right\}^{q^n}
\end{equation}
for any $t\in O$, $\lambda \in D_q$. Setting in (36) $\lambda 
=\log_\delta (z)$ we come to (33). $\blacksquare$

\end{document}